\documentclass[11pt]{article}
\usepackage{amsfonts}
\usepackage{mathrsfs}
\usepackage{amsmath}
\usepackage{amssymb}
\usepackage{graphicx}
\usepackage{epsfig}
\usepackage{epic}
\usepackage{cite}
\usepackage{amsthm,latexsym}
\usepackage{float}

\renewcommand{\paragraph}{\roman{paragraph}}
\newtheorem{theorem}{\scshape \mdseries  Theorem}[section]
\newtheorem{lemma}[theorem]{\scshape \mdseries  Lemma}
\newtheorem{coro}[theorem]{\scshape \mdseries  Corollary}

\newtheorem{prop}[theorem]{\scshape \mdseries  Proposition}
\newtheorem{defi}[theorem]{\scshape \mdseries  Definition}

\newtheorem{example}[theorem]{\scshape \mdseries  Example}

\begin{document}

\title{\sf Determining sets and determining numbers   of finite groups}

\author{Dengyin Wang,\thanks{Corresponding author. E-mail address:3983@cumt.edu.cn. Supported by ``the National Natural Science Foundation of China (No.11571360)".}\ \ \ \ \ \ \ \ Shikun Ou\\
{\small  (School of Mathematics, China University of Mining and Technology, Xuzhou 221116,  China)}\\
 Haipeng Qu\thanks{Supported by ``the National Natural Science Foundation of China
(No.11771258)".}\\
{\small (School of Mathematics and Computer Science, Shanxi Normal University, Linfen 041000, China) } \\ }
\date{}
\maketitle

\noindent {\bf Abstract:}\ \ Let $G$ be a   group. A subset $D$ of $G$ is a determining set of  $G$, if every automorphism of $G$ is uniquely determined by its action on $D$. The determining number of $G$, denoted by $\alpha(G)$, is the  cardinality of a smallest determining set. A generating set of $G$    is a subset such that every element of $G$ can be expressed as the combination, under the group operation, of finitely many elements of the subset and their inverses. The   cardinality of a smallest generating set of $G$, denoted by $\gamma(G)$, is called the  generating number of $G$. A group $G$ is called a DEG-group if $\alpha(G)=\gamma(G)$.

 The main results of this article are as follows. Finite groups  with determining number $0$ or $1$ are classified;  Finite simple groups and   finite nilpotent groups are proved to be DEG-groups;  A finite group is a normal subgroup of a DEG-group and there is an injective mapping from the set all  finite groups   to  the set of finite DEG-groups;   Nilpotent groups of order $n$ which have the maximum determining number are classified;   For any integer $k\geq 2$, there exists a group $G$ such that $\alpha(G)=2$ and $\gamma(G)\geq k$.



\section{Introduction}
\quad\quad When one focusses on the symmetric property of a graph or an algebraic system, it is useful to have in hand a small subset of `key' elements that captures the total symmetric property of the  object. Such a subset will be called a determining set, which originates from the idea of breaking symmetries of the object. The formal definition of determining sets is as follows. Let $G$ be a graph or an algebraic system (such as a linear space, a group, a ring, a Lie algebra, etc.). A subset $D$ of    $G$ is a determining set of $G$, if every automorphism of $G$ is uniquely determined by its action on $D$, i.e., an automorphism  $\sigma$ of $G$ identifies with  another one $\tau$    whenever $\sigma(x)=\tau(x)$ for  any $x\in D$, where an automorphism of a graph is a bijective mapping on the vertex set that preserves adjacency of vertices and an automorphism of an algebraic system is a bijective mapping that preserves all operations in its definition. The determining number of $G$, denoted by $\alpha(G)$,  is the cardinality of a smallest determining set and a determining set of size $\alpha(G)$ is called a minimum determining set.

Determining sets of connected graphs were introduced by Boutin \cite{boutin1}, where ways of finding and verifying determining sets were described. The author also gave natural lower bounds on the determining number of some graphs.   Boutin \cite{boutin2}   studied the determining number of Cartesian products of graphs. Independently, Harary \cite{harary} and Erwin, Harary \cite{erwin} defined an equivalent set and an equivalent number that they called the fixing set and the fixing number, respectively. We refer the reader to \cite{albertson},   \cite{caceres1}, \cite{gibbons} for more works on determining sets of graphs. Determining sets of graphs have application   to distinguishing labeling (see \cite{albertson1}, \cite{albertson2},  \cite{boutin3} or \cite{Imrich}  for details) and   they are frequently used to identify the automorphism groups of   graphs (see \cite{bailey} or \cite{boutin1}).

Determining number of a graph has close connection with another   well-known parameter: the metric dimension or location number of the graph.  For an ordered subset $W =\{w_1, w_2,\ldots, w_k\} $   and a vertex $v$ of $ V(G)$,  the $k$-vector $$r(v|W)=(d(v,w_1), d(v,w_2),\ldots, d(v,w_k))$$ is called the (metric) representation of $v$ with respect to $W$, where $d(v,w_i)$ is the distance of $v$ and $w_i$. The set $W$ is called a resolving set of $G$ if $r(u|W)=r(v|W)$ implies that $u = v$ for all $u, v \in  V(G)$.   A resolving set $S$ of minimum cardinality is a metric basis, and $|S|$, the cardinality of $S$,  is the metric dimension of $G$, which is denoted by $\beta(G)$.  Metric basis and metric dimension have   been widely studied (see  \cite{caceres},\cite{chartrand},   or \cite{bailey} for   review on this parameter),  because of their wide applications to network discovery and verification, robot navigation,  and strategies for mastermind game.
  It has been shown   (see \cite{boutin1},\ \cite{caceres1},\  \cite{erwin}, \cite{harary}) that     $\alpha(G)\leq \beta(G)$ for a connected graph $G$.

  Now, we turn to determining numbers of algebraic systems. For a linear space $\mathbb{V}$ with at least three vectors, a minimum determining set of $\mathbb{V}$  is just a base of $\mathbb{V}$ and the determining number of $\mathbb{V}$  is precisely  the dimension of $\mathbb{V}$.
  In this paper, we  are interested in determining sets and determining numbers of groups.
  \begin{defi} A subset $D$ of a group $G$ is a determining set if  every automorphism of $G$ is uniquely determined by its action on $D$, or equivalently, only the identity automorphism can fixes every element of $D$. The determining number of $G$, denoted by $\alpha(G)$,  is the cardinality of a smallest determining set and a determining set of size $\alpha(G)$ is called a minimum determining set. If $G$ has only the identity automorphism we put $\alpha(G)=0$.\end{defi}

See  $\mathbb{Z}_p^n$ for an example, where   $p$ is a prime number, $\mathbb{Z}_p$ is the quotient group of $\mathbb{Z}$ (the additive group of all integers) to its normal subgroup $p\mathbb{Z}$, and $\mathbb{Z}_p^n$ is the direct product of $n$ copies of $\mathbb{Z}_p$.
\begin{example} $\alpha(\mathbb{Z}_p^n)=n$ for $n\geq 2$.\end{example}

 \noindent{\bf Proof.} The set $\{\varepsilon_1, \ldots, \varepsilon_n\}$, where $\varepsilon_k$ has $\overline{1}$ at the $k$th position and $\overline 0$ elsewhere, is a determining set of  $\mathbb{Z}_p^n$ and thus $\alpha(\mathbb{Z}_p^n)\leq n$. Let $\mathbb{V}$ be the   linear space   over $\mathbb{Z}_p$, which has $\mathbb{Z}_p^n$ as   additive group and the scalar multiple of an element   $\overline{a}$ of  $\mathbb{Z}_p$ with a  vector $(\overline{b_1},\ldots, \overline{b_n})\in \mathbb{V}$ is defined naturally: $$\overline{a}  (\overline{b_1},\ldots, \overline{b_n})=(\overline{ab_1},\ldots, \overline{ab_n}).$$ Since an invertible linear transformation on $\mathbb{V}$ is an automorphism of $\mathbb{Z}_p^n$,     a subset $D$ of  $\mathbb{V}$ is a determining set of $\mathbb{V}$ (as a linear space) if it is a determining set of $\mathbb{Z}_p^n$ (as a group). Thus  $\alpha(\mathbb{Z}_p^n)\geq \alpha(\mathbb{V})=n$, leading to $\alpha(\mathbb{Z}_p^n)=n$. \hfill$\square$

For a permutation group $G$, a well studied parameter, called the base size of $G$ (see \cite{bailey} for more related references) has close relation with determining number of a graph.
\begin{defi} A base for a permutation group $G$, acting faithfully on a finite set $\Omega$, is a subset of $\Omega$ chosen so that its pointwise stabiliser in $G$ is trivial. The base size of $G$ in its action on $\Omega$ is the cardinality of the smallest base for $G$ in this action.\end{defi}

Indeed, the determining number of a graph $ \Gamma$  is   the base size of the automorphism group of $\Gamma$.
  However, determining set and determining number of a group have no  relation with base and base size of a permutation group. In fact, determining set  and determining number of a group can be defined for an arbitrary group, however, base and base size only apply to  permutation groups.
  Even   for a permutation group $G$, there is no relation between a determining set and a base of $G$, since a determining set is a subset of $G$ and a base is a subset of the set on which $G$ acts.

  The problem of  determining all automorphisms of a group is a classic research problem.  When one aims to solve such a problem for a given group $G$, it is very helpful to has a  determining set of $G$ in hand.
  The idea of using  determining sets to study the automorphisms of a group is as follows (see \cite{dyer}, \cite{hillar}, \cite{hua},\cite{san}, \cite{wang} for details): Let $D$ be a   determining set of a group $G$ and  $\sigma$  an arbitrary automorphism of   $G$, if one can find an element $g$  of $G$ such that   ${\rm Inn}(g)\cdot\sigma$    fixes every element of  $D$, then $G$ is proved to have inner automorphisms only, where ${\rm Inn}(g)$ denotes the inner automorphism induced by $g$.  On the contrary,  if  ${\rm Inn}(g)\cdot\sigma$  acts nontrivially on  $D$ for any $g\in G$, then $G$ has out-automorphisms. In this case, choose a $g\in G$ such that    ${\rm Inn}(g)\cdot\sigma$ fixes  the elements of  $D$, as many as possible.  Further, if one can construct     an automorphism     acting  as ${\rm Inn}(g)\cdot\sigma$  on $D$, one succeeds to determine the form of the out-automorphism related to $\sigma$. As a result, $\sigma$ is proved to be the composition of ${\rm Inn}(g^{-1})$ and the   out-automorphism.

  This paper is organized as follows. Some elementary results are given in Section 2. The determining number of the direct product of some groups is studied in Section 3; Finite groups with small   determining number are characterized in Section 4.  In Section 5, we prove that finite simple groups  and finite nilpotent groups are DEG-groups;  there is an injective mapping from the set of all  finite groups  to  the set of  finite DEG-groups. In Section 6, it is proved that, for any integer $k\geq 2$, there exists a group $G$ such that $\alpha(G)=2$ and $\gamma(G)\geq k$. At the last section, we conclude the results obtained in this article and we list six problems for further study.

\section{ Some elementary results}
\quad\quad  In this paper, we only study determining numbers of finite groups. If no other statement,   we denote by  $G$  a finite  group  of order at least $2$. The  identity element of $G$ is written as $1$. For a subset $S$ of $G$, we denote by $|S|$ the number of elements in $S$ and by $\langle S\rangle$ the subgroup of $G$ generated by $S$. The centralizer of $S$ in $G$, written as $C(S)$, is defined by  $$C(S)=\{x\in G:xs=sx, \forall s\in S\}.$$ In particular, we call $C(G)$    the center of $G$. If $C(G)=G$, $G$ is called an abelian group. An automorphism of $G$ is a bijective mapping on $G$ that preserves the   multiplication of $G$. The identity automorphism   refers to the automorphism that fixes all elements of $G$.  Such an automorphism is also called a trivial one. All automorphisms of $G$ form a group, written as  ${\rm Aut}\ G$, under the composition of mappings. An element $x$ of $G$ induces an inner automorphism ${\rm Inn}(x)$ in the way ${\rm Inn}(x):y\mapsto xyx^{-1}, \ \forall y\in G$. For a subgroup $H$ of $G$, let $N_G(H)=\{g\in G:gH=Hg\}$ be the normalizer of $H$ in $G$. If $N_G(H)=G$, $H$ is called a normal subgroup of $G$ and we write $H\lhd G$. A group with at least two elements is called a simple group, if it only has the trivial normal subgroups $G$ and $\{1\}$. If $G$ has a   series   of normal subgroups as  $${\displaystyle \{1\}=Z_{0}\lhd  Z_{1}\lhd  \ldots \lhd  Z_{n}=G},$$ then $G$ is called a nilpotent group, where ${\displaystyle Z_{1}=C(G)}$  and $ {\displaystyle Z_{i+1}}$ is the subgroup such that ${\displaystyle Z_{i+1}/Z_{i}=C(G/Z_{i})}$. A group with order a power of a prime number $p$ is said to be a $p$-group.
 A Klein  four-group, written as $K_4$, is the group $\mathbb{Z}_2^2$, the direct product of a pair of $\mathbb{Z}_2$.

 The fundamental theorem for finite abelian groups (see \cite{smith}) will be frequently used in this article.
 \begin{prop} Let $G$ be a finite abelian group of order at least $2$. Then there is a series of positive integers $\{p_1^{s_1}, \ldots, p_m^{s_m}\}$, where $p_1, \ldots, p_m$ are prime numbers (unnecessary distinct) and $s_1, \ldots, s_m$ are positive integers, such that $G\cong \mathbb{Z}_{p_1^{s_1}}\times \ldots \times \mathbb{Z}_{p_m^{s_m}}$. \end{prop}

Some well known results on nilpotent groups (see \cite{smith}) will be applied when we study groups with equal determining number and generating number.

\begin{prop} The following statements are equivalent for finite groups.\\
(a)\ $G$ is a nilpotent group.\\
(b)\ If H is a proper subgroup of $G$, then $H$ is a proper normal subgroup of $N_G(H)$. \\
(c)\ Every Sylow subgroup of $G$ is normal.\\
(d)\ $G$ is the direct product of its Sylow subgroups.\end{prop}

  \begin{prop} (a) Abelian groups and   finite $p$-groups are nilpotent groups. \\
  (b) A group of order the  square  of a prime number is abelian.  \end{prop}

   \begin{prop} Let $G$ be a   nilpotent  group of order $n$.\\
   (a)\  If $m$ is a divisor of $n$, then $G$ has a subgroup of order $m$.\\
   (b) If $H$ is a normal subgroup of order at least $2$, then $H\cap C(G) $ has nonidentity elements.   \end{prop}
  A generating set of  $G$    is a subset such that every element of $G$ can be expressed as the combination, under the group operation, of finitely many elements of the subset and their inverses. The   cardinality of a smallest generating set of $G$, denoted by $\gamma(G)$, is called the  generating number of $G$.
A moment observation leads to a trivial result.
\begin{prop}  Let $G$ be a group.\\
 (i) A generating set of $G$ is a determining set of $G$ and thus $0\leq\alpha(G)\leq \gamma(G)$.\\
(ii) The image  of a determining set under an automorphism is also a  determining  set.
 \end{prop}

The following proposition  indicates that any given element of a minimum determining set can not be generated by the other elements of the set.

\begin{prop}  Let     $H$ be a subgroup of a finite group $G$ that is generated by a minimum determining set of $G$. Then $\alpha(G)=\gamma(H).$  \end{prop}
\noindent{\bf Proof.} Let $D$ be a minimum  determining set with $ \alpha(G)$ elements that generates  $H$. If   $D$ is proved to be a minimum generating set of $H$, then $\gamma(H)=\alpha(G)$ is derived. Indeed, if $D$    is not a minimum generating set of $H$, then a subset $S$ of $H$ with fewer elements, that is $|S|< |D|$,   can  generate $H$. If $\sigma$ is an automorphism of $G$  fixing all elements of $S$, then it fixes all elements of $D$, which forces $\sigma$ to be the identity automorphism of $G$ and $S$ is also a determining set of $G$, a contradiction. \hfill$\square$

Determining sets of a group $G$ have connection with the center of $G$.
\begin{prop}   The centralizer of a determining set $D$ of a group $G$ is the center of $G$, that is $C(D)=C(G)$.  \end{prop}
\noindent{\bf Proof.} Clearly, $C(G)\subseteq C(D)$. Conversely, if $x\in C(D)$, then $xy=yx, \ \forall y\in D$. In other words, the inner automorphism  ${\rm Inn} (x)$   fixes all elements of $D$. As $D$ is a determining set of $G$, ${\rm Inn} (x)$  fixes all elements of $G$, which implies that $xz=zx,\ \forall z\in G$ and thus $x\in C(G)$. Consequently,  $C(D)=C(G)$.  \hfill$\square$

  Proposition 2.8  shows how   $\alpha(G)$ impacts the size of the automorphism group of $G$.
\begin{prop}  Let  $G$  be a  group of order $n$ with determining number $m$.  Then $|{\rm Aut}\ G|\leq \frac{(n-1)!}{(n-m-1)!} $, the equality holds if and only if $G$ is a cyclic group of prime order or a Klein  four-group. \end{prop}
\noindent{\bf Proof.} Given a minimum determining set $D$ of $G$, arrange the elements in $D$ in some order. By Proposition 2.6, $1\notin D$ and the elements in $D$ are pairwise different. Thus $D$ is an ordered $m$-subset of $G^*$ with $m$ different components, where $G^*=G\setminus \{1\}$. Each automorphic image of $D$ must be an $m$-subset of $G^*$ with $m$ different components and the images are different if   the  automorphisms acting on $D$ are different.  There are totaly $\frac{(n-1)!}{(n-m-1)!}$  different ordered $m$-subsets of $G^*$ with $m$ different components. Thus there are at most $\frac{(n-1)!}{(n-m-1)!}$  different automorphisms of $G$.

If $G$ is a cyclic group of prime order, then   $G$ has precisely $n-1$ automorphisms. The determining number of such a group is  $1$ (if $n\geq 3$) or $0$ (if $n=2$), and thus the equality holds.
If $G$ is a Klein four-group, then $n=4, m=2$ and  $\frac{(n-1)!}{(n-m-1)!}=6$. The automorphism group of $\mathbb{Z}_2^2$ is isomorphic to $  GL(2,\mathbb{Z}_2)$ (see \cite{bid} or \cite{han}), the   group of all $2\times 2$ invertible matrices over $\mathbb{Z}_2$, which has order $6$. Thus the equality holds for $\mathbb{Z}_2^2$.

Conversely, let $D$ be a minimum determining set of $G$ with an element $x$. If $|{\rm Aut}\ G|=\frac{(n-1)!}{(n-m-1)!}$, then the images of $x$ under ${\rm Aut}\ G$ go through all nonidentity elements of  $G$. Since an automorphic image of $x$ has the same order as that of $x$,  all nonidentity elements of $G$  must have a same order, say $p$. If $p$ has a proper division $k$ other than $1$, then $x^k$  has order $\frac{p}{k}$, which implies that   $p$ is a prime number. By Sylow theorem, $G$ must be a $p$-group.  Since a $p$-group has nontrivial center,  there is $\sigma\in {\rm Aut} \ G$ such that $\sigma(x)\in C(G)$. Noting that all automorphisms of $G$ stabilize $C(G)$, we further have   $x\in C(G)$ and   $G=C(G)$, i.e., $G$ is an abelian group.
Recalling that every nontrivial element of $G$ has order $p$, we have $G\cong\mathbb{Z}_p^m$ (thanks to Proposition 2.1 and Example 1.2) and $n=p^m$.   By Corollary 4.3 of \cite{bidwell}, we have $|{\rm Aut}\ G|=\prod_{i=1}^m(p^m-p^{m-i})$. Now, by $$\frac{(n-1)!}{(n-m-1)!}= \prod_{i=1}^m(p^m-p^{m-i})$$ we have $m=1,n=p$, or $p=m=2$. Thus $G$ is a cyclic group of prime order or $G\cong \mathbb{Z}_2^2$, a Klein four-group.\hfill$\square$

\section{Determining number of the direct product of some groups}
\quad\quad Boutin \cite{boutin2} used characteristic matrices to study the determining set of   Cartesian products of graphs.
 Let $G= H_1 \times\ldots\times H_m$ be the direct product of $m$ finite groups, where every element of $G$ is written as a column vector with the $k$th entry an element of $H_k$. Motivated by \cite{boutin2} we define the characteristic matrix of an ordered subset of $G$.
 \begin{defi} Let $X=\{g_1, g_2, \ldots, g_t\}$ be an ordered  subset of $G=H_1 \times\ldots\times H_m$. Define the characteristic matrix $M_X$ to be the $m\times t$ matrix whose $(i,j)$-entry is the $i$th entry of $g_j$. Namely, $M_X$ is the block matrix $[g_1,g_2, \ldots, g_t]$ with the $j$th  column to be $g_j$. \end{defi}

\begin{lemma} Let $G = H_1 \times\ldots\times H_m$ be a direct product of $m$ finite groups   and $X$ an ordered subset of $G$. If $X$ is a determining set of $G$, then each row of $M_X$ contains a determining set of the appropriate factor of $G$. Further, $\alpha(G)\geq {\rm max}\{\alpha(H_i)\} $.\end{lemma}

\noindent{\bf Proof.} Suppose that row $i$ of $M_X$ does not contain a determining set of $H_i$. Then there is a nontrivial automorphism $\varphi_i$ of $H_i$ which fixes all the elements in row $i$ of $M_X$. Then the automorphism of $G$ given by
$$\phi((h_1,\ldots, h_m)^T)=(h_1, \ldots, h_{i-1}, \varphi_i(h_i), h_{i+1}, \ldots, h_m)^T$$  is a nontrivial automorphism of $G$ that fixes  all elements in $X$. Thus $X$ is not a determining set of $G$.
The second assertion follows from the first one immediately.\hfill$\square$

 A known result due to Bidwell will be applied in our next result.
 \begin{prop} {\rm (Theorem 2.2, \cite{bidwell})} \ \
 Let $G = H_1 \times\ldots\times H_m$ be a direct product of $m$ finite groups $H_1, \ldots, H_m$, where no pair of the $H_i\  (1\leq i \leq  m)$ have a common direct factor. Then
 ${\rm Aut}\  G\cong\mathcal{ A }$, where $$\mathcal{A}=\{\left(\begin{array}{ccc}\varphi_{11}&\cdots&\varphi_{1m}\\ \vdots&\ddots&\vdots\\ \varphi_{m1}&\cdots&\varphi_{mm}\end{array}\right): \varphi_{ii}\in{\rm Aut}\  H_i,
\ \ \varphi_{ij}\in {\rm Hom}(H_j,C(H_i))\ {\rm for}\ i\not=j \}.$$ \end{prop}

   Proposition 3.3 indicates that every automorphism of $G$ is  induced by a unique matrix $A=(\varphi_{ij})$  and thus it is written as $\phi_A$, the action of  $\phi_A$   on an element $(h_1,\ldots, h_m)^T$ of $G$  is as follows:
$$\left(\begin{array}{ccc}\varphi_{11}&\cdots&\varphi_{1m}\\ \vdots&\ddots&\vdots\\ \varphi_{m1}&\cdots&\varphi_{mm}\end{array}\right)
\left(\begin{array}{c}h_1\\  \vdots\\ h_m\end{array}\right)=\left(\begin{array}{c}  \varphi_{11}(h_1)\cdots\varphi_{1m}(h_m)\\  \vdots\\ \varphi_{m1}(h_1)\cdots\varphi_{mm}(h_m)\end{array}\right).$$

Let $G $ be as   in Proposition 3.3 and $X$ an ordered subset of $G$. Then, an automorphism  $\phi_A$ of $G$, induced by  $A=(\varphi_{ij})$, fixes $X$ pointwisely, if and only if $AM_X=M_X$, and    $X$ is a determining set of $G$ if and only if only the trivial automorphism $\phi_I$ (induced by the identity matrix $I$) can fix $X$ pointwisely.

Next, we consider the condition under which  $\alpha(H_1\times \ldots\times H_m)= {\rm max}\{\alpha(H_i)\}$. The mapping that sends all elements of $H_j$ to    $1\in H_i$ is called a trivial homomorphism from $H_j$ to $H_i$.
\begin{lemma}  Let $G = H_1 \times\ldots\times H_m$ be a direct product of $m$ finite groups $H_1, \ldots, H_m$, where no pair of the $H_i\  (1\leq i \leq  m)$ have a common direct factor.
 If ${\rm Hom}(H_j,C(H_i))$ only contains the trivial homomorphism for all $ i\not=j$,   then   $\alpha(G)= {\rm max}\{\alpha(H_i)\}$. In particular, if every $H_i$ has trivial center, then $\alpha(G)= {\rm max}\{\alpha(H_i)\}$. \end{lemma}
\noindent{\bf Proof.} By Lemma 3.2,  we have $\alpha(G)\geq {\rm max}\{\alpha(H_i)\}$.
 Suppose ${\rm Hom}(H_j,C(H_i))$ contains only the trivial homomorphism for all $ i\not=j$. We need to prove $\alpha(G)\leq {\rm max}\{\alpha(H_i)\}$.
 Assume   $\alpha(H_1)= {\rm max}\{\alpha(H_i)\}=t$ and let $V_1 $ be
a minimum determining set of $H_1$, ordered arbitrarily. For every other $H_i$, we can choose (since $\alpha(H_i)\leq t$) an ordered subset $V_i$ of $t$ elements of $H_i$ which contains a determining set of $H_i$.
Let $g_j$, for $j=1, \ldots, t$,  be the element of $G$  whose $i$th  entry is the $j$th element of $V_i$. Let $X=\{g_1, \ldots, g_t\}$ with characteristic matrix $M_X=[g_1, \ldots, g_t]$.
If we can prove that $X$ is a determining set of $G$, then $\alpha(G)\leq \alpha(H_1)=t$ and we are done.
Let $\phi$ be an automorphism of $G$ fixing all elements of $X$. We need to prove   $\phi$ is the identity automorphism. By Proposition 3.3, there is an $A=(\varphi_{ij})\in \mathcal{A}$ such that $\phi=\phi_A$. The assumption for $\phi$ indicates that $AM_X=M_X$.
Since ${\rm Hom}(H_j,C(H_i))$ contains only the trivial homomorphism,  the matrix $A$ which induces $\phi_A$ is a diagonal matrix in the form $$A={\rm diag}(\varphi_{11}, \ldots, \varphi_{mm}), \ \ {\rm with}\ \ \varphi_{ii}\in {\rm Aut} \ H_i.$$ It follows from $AM_X=M_X$ that $\varphi_{ii}(x)=x$ for any $x\in V_i$ and for any $i$. As $V_i$ is a determining set of $H_i$, $\varphi_{ii}$ is the identity automorphism of $H_i$.
Therefore $\phi_A$ is the identity automorphism of $G$, which completes the proof of the first assertion. The second assertion follows from the first one immediately. \hfill$\square$

\begin{lemma}  Let $G = H_1 \times\ldots\times H_m$ be a direct product of $m$ finite groups $H_1, \ldots, H_m$.
 If  the orders of $H_1, \ldots, H_m$ are pairwise coprime,   then   $\alpha(G)= {\rm max}\{\alpha(H_i)\}$.  \end{lemma}
\noindent{\bf Proof.} Since the orders of $H_1, \ldots, H_m$ are pairwise coprime,  no pair of the $H_i\  (1\leq i \leq  m)$ have a common direct factor. For $i\not=j$, let    $\varphi_{ij}$ be a homomorphism from $H_j$ to the center of $H_i$, then we have $\frac{|H_j|}{|{\rm Ker} \varphi_{ij}|}=|\varphi_{ij}(H_j)|$. If $|\varphi_{ij}(H_j)|\not=1$, it is a nontrivial common divisor of  $|H_j|$ and $|H_i|$, a contradiction. Hence, $\varphi_{ij} $ is a trivial homomorphism from $H_j$ to the center of $H_i$. By Lemma 3.4, the assertion is proved. \hfill$\square$

\begin{lemma}  Let $G = H_1 \times H_2$ be a direct product of  $H_1$ and $ H_2$, where   $H_1$ and $H_2$ have no common direct factor.
 If   $H_1$ has trivial center and  $\alpha(H_1)=\gamma(H_1)$, and $H_2$ only has the trivial automorphism, then  $\alpha(G)= {\rm max}\{\alpha(H_i)\}$. \end{lemma}
\noindent{\bf Proof.}  By Lemma 3.2,  we have $\alpha(G)\geq {\rm max}\{\alpha(H_i)\}$.
  Suppose $$\gamma(H_1)=\alpha(H_1)=t.$$ As $H_2$ only has  the trivial automorphism, $\alpha(H_2)=0$ and thus $\alpha(G)\geq t$.
    Let $Y=\{v_1, \ldots, v_t\}$ be
a minimum generating set of $H_1$ and let $$Z=\{(v_1,1)^T,\ldots, (v_t,1)^T\},$$ where $1$ is the identity element of $H_2$. If we can prove that $Z$ is a determining set of $G$, then $\alpha(G)\leq t$ and thus $\alpha(G)=t= {\rm max}\{\alpha(H_i)\}$. Let $M_Z=\left(\begin{array}{ccc}v_1&\cdots&v_t\\ 1&\cdots&1\end{array}\right)$ be the characteristic matrix of $Z$.
 Let $\phi$ be an automorphism of $G$ fixing all elements of $Z$.  By Proposition 3.3, there is an $A=(\varphi_{ij})\in \mathcal{A}$ such that $\phi=\phi_A$.  We need to prove $\phi_A$ is the identity automorphism of $G$. Since $H_1$ has trivial center, $\varphi_{12}$ is a trivial homomorphism from $H_2$ to the center of $H_1$.  Since $H_2$ only has the trivial automorphism,   $\varphi_{22}$ is the identity automorphism of $H_2$.
 From  $AM_Z=M_Z$ it follows that $\varphi_{11}(v_j)=v_j$ and $\varphi_{21}(v_j)\varphi_{22}(1)=1$ for all $v_j$. From   $\varphi_{11}(v_j)=v_j,\ \forall v_j$, it follows that $\varphi_{11}$ is the trivial automorphism of $H_1$.  By $\varphi_{21}(v_j)\varphi_{22}(1)=1,\ \forall v_j$, we have $\varphi_{21}(v_j)=1, \ \forall v_j$. Thus, $\varphi_{21}$ is the trivial homomorphism from $H_1$ to the center of $H_2$ since $Y$ is a generating set of $H_1$. Hence, $\phi_A$ is a trivial automorphism of $G$ and thus $Z$ is a determining set of $G$. \hfill$\square$

When we characterize those groups with determining number $1$, we need first to study the determining number of $\mathbb{Z}_{p^k}\times \mathbb{Z}_{p^l}$, where $p$ is a prime number. In fact, by applying Theorem 5.3,   we can also obtain $\alpha(\mathbb{Z}_{p^k}\times \mathbb{Z}_{p^l})=2$. Here, we would like to give a direct proof.
 \begin{lemma}    The determining number of  $\mathbb{Z}_{p^k}\times \mathbb{Z}_{p^l}$ is $2$, where $p$ is a prime number and $k, l$ are positive integers. \end{lemma}
\noindent{\bf Proof.} Let $G=\mathbb{Z}_{p^k}\times \mathbb{Z}_{p^l}$ and suppose $k\leq l$. In this lemma, the elements of $G$ are written as row vectors in  form $(\widetilde{a},  \overline{b })$ with $\widetilde{a}\in\mathbb{Z}_{p^k}$ and $\overline{b }\in \mathbb{Z}_{p^l}$. Clearly, $\alpha(G)\leq 2$ since $G$ can be generated by two elements. The mapping  $\sigma$ defined by $$(\widetilde{a},  \overline{b}) \mapsto (\widetilde{a},  \overline{b+ap^{l-k}}),\ \ \forall \ (\widetilde{a},  \overline{b})\in G $$
is a nontrivial automorphism of $G$, thus  $\alpha(G)\geq 1$. To prove $\alpha(G)=2$ it suffices to prove $\alpha(G)\not=1$.
Suppose, for a contradiction, that $\alpha(G)=1$ and $ (\widetilde{a},  \overline{b})$ is a determining set of $G$. 

 If the order of $\overline{b}$  in  $\mathbb{Z}_{p^l}$  is not smaller than that of  $\widetilde{a}$  in  $\mathbb{Z}_{p^k}$, then  $ (\widetilde{a},  \overline{b})$ has the same order as that of $\overline{b}$. The image of  $ (\widetilde{a},  \overline{b}) $ under $\sigma$ has the same order as that of $ (\widetilde{a},  \overline{b}) $. Thus  $(\widetilde{a},  \overline{b+ap^{l-k}})$ has the same order as that of  $\overline{b}$, which implies that $\overline{b+ap^{l-k}}$ and $\overline{b}$ have the same order. It is easy to see that
the elements in  $\mathbb{Z}_{p^l}$  with a same order are in a same orbit under the action of automorphisms of $\mathbb{Z}_{p^l}$. There is an automorphism $\tau$ of $\mathbb{Z}_{p^l}$ such that $\tau(\overline{b+ap^{l-k}})= \overline{b}$. Let $\tau'$ be the automorphism on $G$ sending any $(\widetilde{a},  \overline{b })$ to $(\widetilde{a},  \tau(\overline{b})) $. Then the nontrivial automorphism $\tau'\cdot\sigma$ of $G$, sending  $(\widetilde{1},  \overline{0})$ to $(\widetilde{1},  \tau(\overline{p^{l-k}})) $, fixes  $(\widetilde{a},  \overline{b }) $, which is a contradiction.

 The condition $l\geq k$ allows us to define an automorphism $\lambda$ of $G$ as: $$\lambda:\  (\widetilde{a},  \overline{b}) \mapsto (\widetilde{a+b},  \overline{b}),\ \ \forall \ (\widetilde{a},  \overline{b})\in G.$$
  If the order of $\widetilde{a}$  in  $\mathbb{Z}_{p^k}$ is larger than that of $\overline{b}$  in  $\mathbb{Z}_{p^l}$, then, as the discussion at the above paragraph, there is an automorphism $\pi$ of   $\mathbb{Z}_{p^k}$ such that  $\pi(\widetilde{a+b})=\widetilde{a}$. Let $\pi'$ be the automorphism on $G$ sending any $(\widetilde{a},  \overline{b })$ to $(\pi(\widetilde{a}),   \overline{b}) $. Then the nontrivial automorphism $\pi'\cdot\lambda$ of $G$, sending $(\widetilde{0},  \overline{1 })$ to $(\pi(\widetilde{1}),   \overline{1}) $, fixes  $(\widetilde{a},  \overline{b }) $, which is a contradiction.\hfill$\square$

\section{Finite groups with small determining numbers}
\quad\quad As the authors of \cite{caceres1} pointed out,   a graph $\Gamma$  has determining number $0$ if and only if $\Gamma$ only has the identity automorphism, i.e., $\Gamma$ is an asymmetric graph. However,   almost all graphs are asymmetric (see Corollary 2.3.3, \cite{godsil}), hence  have determining number $0$. It is still an open problem to give a graphic characterization for graphs $\Gamma$ with  $\alpha(\Gamma)=0$.
If one turns to consider the corresponding problem for a finite group, the situation is quite different, only one group has determining number $0$.
\begin{theorem}  For a   finite group $G$, $\alpha(G)=0$ if and only if $G$ is a cyclic group of order $2$.  \end{theorem}
\noindent{\bf Proof.} If  $G$ is a cyclic group of order $2$, then $G$ only has the trivial automorphism and thus $\alpha(G)=0$.
Conversely, if $\alpha(G)=0$, then $G$ is an abelian group since the inner automorphism induced by a noncentral element of a non-abelian group is a nontrivial automorphism of the group.
When $G$ is abelian, the mapping $\sigma$ on $G$ defined by $\sigma(g)=g^{-1},\  \forall g\in G$, is an automorphism of $G$. From $\alpha(G)=0$ it follows that $g= g^{-1}$ for all $g\in G$, i.e., every nonidentity element of $G$ has order $2$.
Thus $G$ is an abelian $2$-group. The fundamental theorem for finite abelian groups applied to $G$ implies that  $G\cong\mathbb{Z}_2^m$ for some $m$.  Example 1.2 indicates that $\alpha(\mathbb{Z}_2^m)=m$ for $m\geq 2$, from which it follows that $m=1$ and thus $G\cong\mathbb{Z}_2$, as required.  \hfill$\square$

Theorem 4.1 gives a group,  that is $\mathbb{Z}_2$, whose determining number is not equal to its generating number, although the example seems somewhat trivial.

An algebraic characterization of those graphs with $\alpha(\Gamma) = 1$ was obtained by Erwin and Harary \cite{erwin} as follows: Let $\Gamma$ be a nonidentity graph. Then $\alpha(\Gamma) = 1$  if and only if $\Gamma$ has an orbit of cardinality $|{\rm Aut} \ \Gamma|$. However, the problem of  giving a graphic  characterization for graphs
$\Gamma$ with $\alpha(\Gamma) = 1$ is still open.
If we turn to groups,  the corresponding problem is not difficult.

\begin{lemma}  If a   finite group $G$ is not an abelian group, then $\alpha(G)\geq 2$.  \end{lemma}
\noindent{\bf Proof.} By Theorem 4.1,   $\alpha(G)\geq 1$. We only need to prove that $\alpha(G)\not=1$. Suppose for a contradiction that $\alpha(G)=1$ and $\{x\}$ is a determining set of $G$.

If $x\in C(G)$, then the inner automorphism induced by a noncentral element of $G$ fixes $x$, however, the automorphism is not the identity automorphism, which is a contradiction.

If $x\notin C(G)$, then the  inner automorphism induced by $x$ fixes $x$ and it is not  the identity automorphism,  leading to a contradiction.\hfill$\square$

\begin{theorem} A  finite group $G$ has determining number $1$   if and only if $G$ is a cyclic group of order at least  $3$.  \end{theorem}

\noindent{\bf Proof.} The sufficient direction is obvious. Let $G$ be a finite group with $\alpha(G)=1$. By Lemma 4.2, $G$ is abelian. The fundamental theorem for finite abelian groups applied to $G$ gives that
$G\cong \mathbb{Z}_{p_1^{s_1}}\times\ldots \times  \mathbb{Z}_{p_m^{s_m}}$, where $p_1, \ldots, p_m$ are prime numbers (unnecessary different) and $s_1, \ldots, s_m$ are  positive integers.
If a pair of prime numbers in $\{p_1, \ldots, p_m\}$, say $p_1, p_2$, are equal, then
$G$ has  $\mathbb{Z}_{p_1^{s_1}}\times  \mathbb{Z}_{p_1^{s_2}}$ as a direct factor. By Lemma 3.2 and Lemma 3.7, we have $\alpha(G)\geq \alpha(\mathbb{Z}_{p_1^{s_1}}\times  \mathbb{Z}_{p_1^{s_2}})=2$, which is a contradiction. Thus   $p_1, \ldots, p_m$ are pairwise different and therefore $G$ is a cyclic group of order $ p_1^{s_1}\cdots p_m^{s_m}$, which completes the proof.\hfill$\square$

By  Theorem 4.1 and Theorem 4.3, a corollary follows  immediately, which gives a family of groups with equal determining number and generating number. In the next section, more families of such groups will be given.
\begin{coro} If  the generating number of a finite group is $2$, then the determining number of the group is $2$.   \end{coro}

 \section{Finite groups whose determining number and generating number are equal}
 \quad\quad In this section, we study those  groups  with  equal determining number and generating number. For convenience,  a group $G$ will be called a DEG-group (shorthand for `determining number equals generating number') if $\alpha(G)=\gamma(G)$. It should be pointed out that a minimum generating set in a DEG-group $G$  must  be a minimum determining set, however, a minimum determining set sometimes fails to be a minimum generating set.
 See the cyclic group $\mathbb{Z}_6$ for an example.   It is easy to see that $\alpha(\mathbb{Z}_6)=\gamma(\mathbb{Z}_6)=1$   and $\{\overline 2\}$ is a minimum determining set, however,  it is not a    generating set.

By the above section, we have known that all groups $G$ with $\gamma(G)\leq 2$, except for $\mathbb{Z}_2$, are DEG-groups. Now, we try to find more such groups.

 \begin{theorem} A finite simple group  of order at least $3$  is a DEG-group. \end{theorem}
\noindent{\bf Proof.} Let $G$ be     a finite simple group. Then  $\gamma(G)\leq 2$, which is proved by Martino et.al. in  \cite{mart} (or see \cite{barl}, page 195-235).
If $G$ is a cyclic group of order at least $3$, then it follows from Theorem 4.3 that $\alpha(G)=\gamma(G)=1$. If $\gamma(G)=2$, then   Corollary 4.4 confirms the  assertion. \hfill$\square$

Before  the study on an arbitrary nilpotent group, we focus on two   families of special nilpotent groups, $p$-groups and abelian groups.

\begin{lemma} Let $G$ be a finite group, $p$   a prime number and $M$ a normal subgroup of $G$. If the index of $M$ in $G$ is $p$ and $M\cap C(G)$ contains an element of order $p$, then $G$ has a nontrivial automorphism which
fixes all elements of $M$. \end{lemma}
\noindent{\bf Proof.} Suppose $z\in M\cap C(G)$ is of order $p$. Let $a\notin M$ and $K$ be the subgroup generated by $a$. Then $G=KM$ and every element $g$ of $G$ can be written as $g=a^ix$ with $0\leq i\leq p-1$ and $x\in M$.
Define a mapping $\sigma$ from $G$ to itself by  $$\sigma:\ \ a^ix\mapsto (az)^ix,\ \    0\leq i\leq p-1,\ x\in M.$$
It is easy to see that $\sigma$ is a bijection,  it fixes all elements of $M$ and it is not the identity mapping. To complete the proof, we only need to prove that $\sigma$ preserves the   operation of $G$.
For any $g_1, g_2$ of $G$, write them as $$g_1=a^ix, \ \ g_2=a^jy, \ \ {\rm where}\  0\leq i,j\leq p-1,\ \  x,y\in M.$$ Then $g_1g_2=a^{i+j}x'y$ with $x'=a^{-j}xa^j\in M$. If $i+j\leq p-1$, then  by the definition of $\sigma$ one can easily find that
$$\sigma(g_1g_2)=(az)^{i+j}x'y=\sigma(g_1)\sigma(g_2).$$
If $i+j\geq p$, then $g_1g_2=a^{i+j-p}(a^px'y)$, where $0\leq i+j-p\leq p-1$, $a^px'y\in M$. By the definition of $\sigma$, we also have $$\sigma(g_1g_2)=(az)^{i+j-p}(a^px'y)=\sigma(g_1)\sigma(g_2).$$ \hfill$\square$

 \begin{theorem} A finite $p$-group  of order at least $3$  is a DEG-group. \end{theorem}
\noindent{\bf Proof.} Let $G$ be the $p$-group and let  $D$ be a minimum determining set of $G$. If we can prove that $D$ is a generating set of $G$, then $\alpha(G)\geq \gamma(G)$ and we are done. Suppose, for a contradiction, that  $D$ generates a proper subgroup $H$ of $G$. Let $M$ be a maximal subgroup of $G$ that contains $H$. By (b) of Proposition 2.2, $M$ is a proper subgroup of $N_G(M)$. As $M$ is a maximal subgroup of $G$, $N_G(M)=G$, which  implies that $M$ is a normal subgroup of $G$ and the index of $M$ in $G$ is $p$.  Further, by (b) of Proposition 2.4,  $|M\cap C(G)|\geq 2$ and thus there is an element $z\in M\cap C(G)$ which has order $p$. By Lemma 5.2, $G$ has a nontrivial automorphism fixing all elements of $M$, which contradicts the condition that $D$ is a determining set of $G$.  \hfill$\square$

\begin{lemma} Let $G$ be a finite abelian group, $p$   an odd prime and $M$ a  subgroup of $G$. If  the index of $M$ in $G$ is $p$, then $G$ has a nontrivial automorphism that
fixes all elements of $M$. \end{lemma}
\noindent{\bf Proof.}  Let $P$ be the Sylow $p$-subgroup of $G$. If $|P\cap M|\not=1$, then $M$ has an element of order $p$ and the proof is established by Lemma 5.2. If  $|P\cap M|=1$ then $P$ is a cyclic group of order $p$ and
$G$ is the direct product of $P$ and $M$. Thus every element $g$ of $G$ can be uniquely written as $g=xz$ with $x\in P$ and $z\in M$. Let $\sigma$    be the mapping on $G$ defined by  $$xz\mapsto x^{-1}z  \ \ {\rm where}\    x\in P, \ z \in M.$$
It is easy to see that $\sigma$ is an automorphism of $G$ that fixes all elements of $M$. Since $p\not=2$, there is $x\in P$ such that $x\not=x^{-1}$, which implies  $\sigma$ is nontrivial.
  \hfill$\square$
\begin{theorem} A finite abelian group  of order at least $3$  is a DEG-group. \end{theorem}
\noindent{\bf Proof.} Let $G$ be the abelian group. If we can find a minimum determining set which  generates  $G$, then $\alpha(G)\geq \gamma(G)$ and we are done. Let $D$ be a minimum determining set such that $|\langle D\rangle|\geq |\langle S\rangle|$ for any minimum determining set $S$ of $G$.  We aim to prove that $\langle D\rangle=G$. Suppose, for a contradiction, that  $D$ generates a proper subgroup $H$ of $G$. Let $M$ be a maximal subgroup of $G$ that contains $H$.
Then the index of $M$ in $G$ is a prime number, which is denoted by $p$.

 If $p$ is odd, then it follows from Lemma 5.4 that $G$ has     a nontrivial automorphism that
fixes all elements of $M$, which is a contradiction.

For the case when $p=2$, if the order of $M$ is even, then $M$ has an element of order $2$ (by (a) of Proposition 2.4) and thus we can derive a contradiction by applying Lemma 5.2. Suppose $p=2$ and the order of $M$ is odd. Then $G$ is the direct product of a two-elements group, say $\{1, x\}$, and $M$. Let $D'=\{xd:d\in D\}$. Note that $x, d$ are both powers of $xd$ for $d\in D$. Hence, $\langle D\rangle\subsetneq \langle D'\rangle$, and  if an automorphism of $G$ fixes $xd$, then it fixes both $x$ and $d$, which further implies that $D'$ is also a minimum determining set of $G$, a contradiction to the choice of $D$. \hfill$\square$

\begin{lemma} Let $G=K\times H$ be a direct product of two subgroups of a finite group $G$, where the orders of $K$ and $H$ are coprime. If    $\gamma(K)\geq \gamma(H)$ and $K$ is a DEG-group, then $G$ is a DEG-group. In particular, $G$ is a DEG-group if so are $K$ and $H$. \end{lemma}
\noindent{\bf Proof.} By Lemma 3.5, $\alpha(G)={\rm max}\{\alpha(K), \alpha(H)\}=\alpha(K)$. Suppose $\gamma(K)=k\geq \gamma(H)=l$. We try to prove $\gamma(G)=\gamma(K)=k$. Obviously,   $\gamma(G)\geq k$. Let $\{x_1,\ldots, x_k\}$, $\{y_1, \ldots, y_l\}$ respectively be  a generating set of $K$ and $H$. Since the order of $x_i$ and $y_j$ are coprime and they are commutative, $x_i$ and $y_j$ are both some powers of $x_iy_j$. Then it is easy to see that $\{x_1y_1, \ldots, x_ly_l, x_{l+1}, \ldots, x_k\}$ is a generating set of $G$  and hence $\gamma(G)\leq k$. Further, we have $\gamma(G)=k$.  Thus $$\gamma(G)=k=\gamma(K)=\alpha(K)=\alpha(G),$$ showing that $G$ is a DEG-group. The second assertion is obvious.\hfill$\square$

With Theorem 5.3 and Lemma 5.6 in hand, we are now ready to study an arbitrary finite nilpotent group.

\begin{theorem} A finite nilpotent group  of order at least $3$  is a DEG-group. \end{theorem}
\noindent{\bf Proof.} Let $G$ be the nilpotent group. By Proposition 2.2, $G$ can be written as the direct product of its Sylow subgroups:
$$G=  H_1 \times\ldots\times H_m.$$ The orders of these Sylow subgroups $H_i$ are pairwise coprime. Arrange $H_i$ such that  $\gamma(H_1)\geq\gamma(H_2)\geq \ldots \geq \gamma(H_m)$.
We proceed by induction on $m$ to prove the result. If $m=1$, then $G=H_1$ is a $p$-group and the result is proved by Theorem 5.3. Assume the result holds for $m-1$. Let $G_1= H_1 \times\ldots\times H_{m-1}$. The induction hypothesis says that $G_1$ is a DEG-group. Now, $G=G_1\times H_m$ and $\gamma(G_1)=\alpha(G_1)=\gamma(H_1)\geq \gamma(H_m)$. By Lemma 5.6, $G$ is a DEG-group.   \hfill$\square$

Next, we consider the  proportion of finite DEG-groups in all finite groups. For a finite group $H$ with $\gamma(H)=l$, let $Q_H=\mathbb{Z}_p^l\times H$, where $p$ is the minimum prime number satisfying  $p>|H|$. Then $Q_H$   is uniquely defined by $H$, which  will be called    a  tight cover of $H$.

\begin{lemma} For a finite group $H$,   there is a DEG-group $G$ that contains $H$ as a normal subgroup.  \end{lemma}
\noindent{\bf Proof.}  Suppose $\gamma(H)=l$ and let $Q_H$ be the tight cover of $H$. Then $Q_H$ contains $H$ as a normal subgroup. Since the order of $\mathbb{Z}_p^l$ and $H$ are coprime, $\gamma(\mathbb{Z}_p^l)=\gamma(H)=l$ and $\mathbb{Z}_p^l$ is a DEG-group, by Lemma 5.6, $Q_H$ is a DEG-group.   \hfill$\square$

 \begin{theorem} There is an injective mapping from the set of all  finite groups  to the set of finite DEG-groups. \end{theorem}
\noindent{\bf Proof.} View  isomorphic groups   as a same group. Let $\Sigma$ be the set of finite groups and $\Psi$ the set of finite DEG-groups.
Let $\eta$ be the mapping from $\Sigma$ to $\Psi$ sending any  $X\in\Sigma $ to $Q_X\in\Psi$. To complete the proof, we need to prove that $\eta$ is injective. It suffices to prove that $Q_{X_1}\cong Q_{X_2}$ implies $X_1\cong X_2$ for $X_1, X_2\in \Sigma$.  Let $p_i$ be the minimum prime number larger than $|X_i|$ and $\gamma(X_i)=l_i$. If  $p_1\not= p_2$, say $p_1>p_2$,  since $p_1$ is a divisor of $|Q_{X_2}|$ we have $p_1$ is a divisor of   $p_2^{l_2}$, a contradiction. Hence, $p_1=p_2$. Similar discussion leads to $l_1=l_2$.
Finally, we have  $$X_1\cong Q_{X_1}/\mathbb{Z}_{p_1}^{l_1}\cong Q_{X_2}/\mathbb{Z}_{p_2}^{l_2}\cong X_2,$$ as required. \hfill$\square$

Next, we consider the relation between $\alpha(G)$ and the  order of  $G$.
Denote by $\chi(G)$ the number of prime divisors (with multiplicity) of $|G|$. More definitely, $\chi(G)=  m_1+\ldots+m_2$ if $$|G|=p_1^{m_1}p_2^{m_2}\cdots p_s^{m_s}$$ is the decomposition of $|G|$ into the product of some  prime divisors $p_1, \ldots, p_s$.

\begin{theorem} Let $G$ be a  finite group.   Then  $\alpha(G)\leq \chi(G)$. If $\alpha(G)=\chi(G)$ then all the following assertions hold  for $G$.\\
 (i) For any subgroup $H$ of $G$, $\gamma(H)=\chi(H)$.\\
 (ii)  Every nonidentity element of $G$  has a prime order.\\
 (iii)   $G$ is a nilpotent   group   if and only  if $G$ is the direct product of $\chi(G)$ copies of cyclic group of prime order.
 \end{theorem}
\noindent{\bf Proof.} Since $\alpha(G)\leq \gamma(G)$, if we can  prove  $\gamma(G)\leq \chi(G)$, then $\alpha(G)\leq \chi(G)$ is derived. We proceed by induction on $\chi(G)$ to achieve the goal. If $\chi(G)=1$, then $G$ is a cyclic group and thus $ \gamma(G)\leq1$, as required.
Suppose the assertion holds for groups whose orders have at most $\chi(G)-1$ prime divisors. For group $G$, let $K$ be a  maximal subgroup of $G$, then the order of $K$ has at most $\chi(G)-1$ prime divisors, thus, by the induction hypothesis, $\gamma(K)\leq \chi(K)$.  Since the union of  a generating set of $K$ together with an element not in $K$ generate  $G$, we have  $$\gamma(G)\leq \gamma(K)+1\leq \chi(K)+1\leq \chi(G),$$ as required.

For the rest of the proof, assume $\alpha(G)=\chi(G)$.  Let $\Sigma$ be the set of  subgroups $H$ of $G$  with $\gamma(H)\not=\chi(H)$. If $\Sigma=\emptyset$, then (i) is proved.   Suppose for a contradiction that $\Sigma\not=\emptyset$ and let $H$ be a subgroup in $\Sigma$ with the maximum order. The condition $\alpha(G)= \chi(G)$ implies that $H$ is a proper subgroup of $G$. Let $x$ be an element of $G$ outside $H$ and let $K=\langle H, x\rangle$.  Then $\gamma(K)\leq \gamma(H)+1$ and $\chi(K)\geq \chi(H)+1$, from which it follows  $\chi(K)\geq \chi(H)+1>\gamma(H)+1\geq\gamma(K)$. Thus $K\in \Sigma$, which is a contradiction to the choice of $H$. The proof for (i) is completed.

Let $z$ be a nonidentity element of $G$ and let $Z$ be the cyclic group generated by $z$. By (i), $\chi(Z)=\gamma(Z)=1$ and thus $Z$ has  a prime order, which completes the proof for (ii).

The  sufficiency for (iii) is obvious.

For the necessity,  we first consider the case when $G$ is abelian.
Let $G$ be a finite abelian  group with $\alpha(G)=\chi(G)$. Then all prime divisors of $|G|$ must be identical, say $p$. To see this, we remind that   $x_1x_2\in G$ has order $p_1p_2$ if the order of $x_i$ is $p_i$ and   $p_1, p_2$ are   distinct prime numbers. Applying   Proposition 2.1 to $G$, we immediately have $G\cong\mathbb{Z}_p^k$ with $k=\chi(G)$.

Let $G$ be an arbitrary finite nilpotent group  with $\alpha(G)=\chi(G)$. (d) of  Proposition 2.2 indicates that $G$ can be decomposed into the direct product of the Sylow subgroups of $G$. By (ii) of this theorem, we find that there is only one direct product factor  in the decomposition of $G$ (otherwise, the order of the product of two elements of prime orders from distinct factors is not a prime number).
Thus $G$ is a $p$-group with $p$ a prime. If we can prove that $G$ is abelian, then we are done. For any nonidentity elements $x, y\in G$, let $H=\langle x,y\rangle$. By (i) of this theorem, $\chi(H)=\gamma(H)\leq 2$. If $\chi(H)=1$, then $H$ is a cyclic group of order $p$ and thus $xy=yx$. If $\chi(H)=2$, then $H$ has order $p^2$ and thus $xy=yx$ (thanks to (b) of Proposition 2.3). \hfill$\square$

\section{Groups that are not DEG}
\quad\quad The cyclic group of order $2$ is a trivial example which is not DEG. In this section, we try to find some other such groups.
 Our main result in this section shows that the difference of $\gamma(G)$ and $\alpha(G)$ can be much large. Before giving such groups, we introduce a known result about the automorphisms of a solvable group. Let $F$ be a field of characteristic not $2$ and $n\geq 2$. Let $F^*$ be the multiplicative group of $F$ and $T_n^*(F)$ be the multiplicative group of all invertible upper triangular $n\times n$ matrices over $F$. Denote   by $E_{i,j}$ the square matrix   whose $(i,j)$-entry is $1$ and all other entries are $0$. The $n\times n$ identity matrix   is written as $I_n$.

 \begin{prop}{\rm (Theorem 1, \cite{zhang})} A mapping $f:T_n^*(F)\to T_n^*(F)$ is a group automorphism if and only if there is a matrix $Q\in T_n^*(F)$ such that either

 (i)\ \  $f(A)=\psi(A)QA^\sigma Q^{-1}, $ \ \  $\forall\ A\in  T_n^*(F),$  \ \ \ or

 (ii)\ \  $f(A)=\psi(A)Q[J(A^\sigma)^{-T}J] Q^{-1}, $ \ \  $\forall\ A\in  T_n^*(F),$\\
 where $\sigma$ is a field automorphism of $F$, $A^\sigma=(\sigma(a_{ij}))$ for $A=(a_{ij})$, $A^{-T}$ is the transpose inverse of $A$, $J=\sum_{i=1}^n E_{i,n-i+1}$, $\psi:T_n^*(F)\to F^*$ is a homomorphism that satisfies
 $$\{\psi(xI_n)\sigma(x):x\in F^*\}=F^*\ \ \ \ {\rm and}\ \  \ \ \{x:\psi(xI_n)\sigma(x)=1\}=\{1\}.$$
 \end{prop}

 Let $p$ be an odd prime number, $\mathbb{F}_p$ be a finite field of $p$ elements and let $ST_{p-2}(\mathbb{F}_p)$ be    the  subgroup of $T_{p-2}^*(\mathbb{F}_p)$ of matrices with determinant $1$.
 Applying Proposition 6.1, one can easily characterize the automorphisms of $ST_{p-2}(\mathbb{F}_p)$.

\begin{lemma}  A mapping $f:ST_{p-2}(\mathbb{F}_p)\to ST_{p-2}(\mathbb{F}_p)$ is an   automorphism if and only if there is a matrix $Q\in T_{p-2}^*(\mathbb{F}_p)$ such that either

 (i)\ \  $f(A)=QA Q^{-1}, $ \ \  $\forall\ A\in  ST_{p-2}(\mathbb{F}_p),$  \ \ \ or

 (ii)\ \  $f(A)=Q[J A^{-T}J] Q^{-1}, $ \ \  $\forall\ A\in  ST_{p-2}(\mathbb{F}_p),$\\
 where  $A^{-T}$ is the transpose inverse of $A$, $J=\sum_{i=1}^{p-2} E_{i, p-i-1}$.
 \end{lemma}
\noindent{\bf Proof.} Recall that $\mathbb{F}_p^*$ is a cyclic group of order $p-1$ and $\mathbb{F}_p$ only has the trivial   field automorphism. Thus $a^{p-1}=1$ and $a^{p-2}=a^{-1}$ for any $a\in \mathbb{F}_p^*$. Denote by $d_A$ the   determinant of    $A\in T_{p-2}^*(\mathbb{F}_p)$.
 For any $A\in T_{p-2}^*(\mathbb{F}_p)$, the    determinant of $d_AA$ is $d_A^{p-2}d_A=1$, showing that $d_AA\in ST_{p-2}(\mathbb{F}_p)$. Thus any $A\in T_{p-2}^*(\mathbb{F}_p)$ can be uniquely written as $A=d_A^{-1}(d_AA)$ with $d_AA\in  ST_{p-2}(\mathbb{F}_p)$. Let  $f:ST_{p-2}(\mathbb{F}_p)\to ST_{p-2}(\mathbb{F}_p)$ be an   automorphism  of $ST_{p-2}(\mathbb{F}_p)$.
 Extend it to a mapping $\bar f:   T_{p-2}^*(\mathbb{F}_p)\to   T_{p-2}^*(\mathbb{F}_p)$ as $$\bar f(A)=d_A^{-1}f(d_AA), \ \ \forall    A=d_A^{-1}(d_AA)\in T_{p-2}^*(\mathbb{F}_p).$$
It is easy to see that $\bar f$ is an automorphism of  $T_{p-2}^*(\mathbb{F}_p)$. Then there is a matrix $Q\in T_{p-2}^*(\mathbb{F}_p)$ such that either

 (i)\ \  $\bar f(A)=\psi(A)QA Q^{-1}, $ \ \  $\forall\ A\in   T_{p-2}^*(\mathbb{F}_p),$  \ \ \ or

 (ii)\ \  $\bar f(A)=\psi(A)Q[J A^{-T}J] Q^{-1}, $ \ \  $\forall\ A\in   T_{p-2}^*(\mathbb{F}_p),$\\
 where $\psi:T_{p-2}^*(\mathbb{F}_p)\to \mathbb{F}_p^*$ is a homomorphism that satisfies the condition shown in Proposition 6.1.
If     case (ii) happens for $\bar f$, then  $$ d_A^{-1}f(d_AA)= \psi(A)Q[J A^{-T}J] Q^{-1},\ \  \forall\ A\in   T_{p-2}^*(\mathbb{F}_p).$$
By comparing the determinants we have $\psi(A)=d_A^{-2},  \forall\ A\in   T_{p-2}^*(\mathbb{F}_p).$
Consequently,  $$ f(d_AA)= d_A^{-1}Q[J A^{-T}J] Q^{-1},\ \  \forall\ A\in   T_{p-2}^*(\mathbb{F}_p).$$
Since $d_X=1$  for   $X\in   ST_{p-2}(\mathbb{F}_p)$,  we have   $$ f(X)=Q[J X^{-T}J] Q^{-1}, \ \ \forall\ X\in   ST_{p-2}(\mathbb{F}_p). $$
If   case (i) happens for $\bar f$, a similar discussion leads to $$ f(X)=QXQ^{-1}, \ \ \forall\ X\in   ST_{p-2}(\mathbb{F}_p),$$ where $Q\in T_{p-2}^*(\mathbb{F}_p)$. \hfill$\square$

The main result of this section is as follows.

\begin{theorem} For any integer $k\geq 2$, there exists a group $G$ such that $\alpha(G)=2$ and $\gamma(G)\geq k$.\end{theorem}

\noindent{\bf Proof.} Let $p$ be a prime number with $p\geq k+3$ and let $G=ST_{p-2}(\mathbb{F}_p)$ consisting of all $(p-2)\times (p-2)$ upper triangular matrices over $\mathbb{F}_p$ with determinant $1$.
The mapping $$\tau:ST_{p-2}(\mathbb{F}_p) \to (\mathbb{F}_p^*)^{p-3}$$ defined by
$$\tau(A)=(a_{11},  \ldots, a_{p-3,p-3}), \ \ \forall A=(a_{ij})\in ST_{p-2}(\mathbb{F}_p)$$ is a surjective  homomorphism from  $ST_{p-2}(\mathbb{F}_p)$ to  $(\mathbb{F}_p^*)^{p-3}$. The kernel  of $\tau$, written as $K$,  is   the normal subgroup of $G$ consisting of all $(a_{ij})\in ST_{p-2}(\mathbb{F}_p)$ with $a_{ii}=1, \ \forall \ i$.
 Thus $G/K$ is isomorphic to $(\mathbb{F}_p^*)^{p-3}$, the direct product of $p-3$ copies  of $\mathbb{F}_p^*$. Since $\gamma((\mathbb{F}_p^*)^{p-3})=p-3$, we have
 $$\gamma(G)\geq \gamma(G/K)=\gamma((\mathbb{F}_p^*)^{p-3})=p-3\geq k.$$

  Let $B={\rm diag}(1,2,\ldots, p-2)$ and let $C= I_{p-2}+\sum_{i=1}^{p-3}E_{i,i+1}$. Since $\prod_{k=1}^{p-2}k=1$ in $\mathbb{F}_p^*$, we have $B\in G$. Thus $\Omega=\{B,C\}$ is a  subset of $G$ with two elements.
 To complete the proof we only need to prove that $\Omega$ is a determining set of $G$.

 Suppose an automorphism $\phi$ of $G$ fixes $B$ and $C$, respectively.  By Lemma 6.3,  there exists a matrix  $Q\in T_{p-2}^*(\mathbb{F}_p)$ such that either

 (i)\ \  $\phi(A)=QA Q^{-1}, $ \ \  $\forall\ A\in  G,$  \ \ \ or

 (ii)\ \  $\phi(A)=Q[J A^{-T}J] Q^{-1}, $ \ \  $\forall\ A\in  G.$

 We claim that case (ii) is impossible. Otherwise, $\phi(B)=Q[J B^{-T}J] Q^{-1}$ is a matrix whose $(1,1)$-entry is  $(p-2)^{-1}$ and thus $(p-2)^{-1}=1$, absurd.
 Consequently, there is  matrix  $Q\in T_{p-2}^*(\mathbb{F}_p)$ such that   $$\phi(A)=QA Q^{-1},  \ \   \forall\ A\in  G.$$ From $\phi(B)=B$ it follows that $QB=BQ$, which implies that
 $Q$ must be a diagonal matrix (noting that the diagonal entries of $B$ are pairwise distinct). Further, it follows from $\phi(C)=C$ that $QC=CQ$, which implies that $Q$ is a scalar matrix. Hence $\phi(A)=A$ for  $\forall\ A\in  G,$ i.e., $\phi$ is a trivial automorphism and $\alpha(G)\leq 2$. Since $G$ is not abelian, $\alpha(G)\geq 2$.  Therefore, $\alpha(G)=2$. \hfill$\square$

\section{ Conclusions and further research problems}

 \quad\quad In this article, finite groups $G$ with determining number $0$ or $1$ are classified. It has been shown that many  finite groups have determining number $2$ because finite groups with $\gamma(G)=2$ are such groups, including all non-cyclic finite simple groups. A further   research problem is:\\
{\bf Problem 1:}{\it \  Classify the finite groups whose determining number are  $2$}.

 Theorem 5.9 indicates    there is an injective mapping from the set of all finite groups  to  the set of   finite DEG-groups. Here, a more interesting problem is:\\
 {\bf Problem 2:}\ {\it  In the set of groups of  order $n$ (isomorphic groups are viewed as the same), what is the value of $$  \lim_{n\rightarrow \infty}\frac{|\mathbb{D}_n|}{|\mathbb{T}_n|},$$  where $\mathbb{D}_n$ is the set of DEG-groups of order $n$  and $\mathbb{T}_n$ is the set of all groups of order $n$.}

 Theorem 5.10 has characterized finite nilpotent groups which satisfy $\alpha(G)=\chi(G)$. A general problem is:\\
    {\bf Problem 3:}\ {\it Characterize  an arbitrary finite   group $G$ which satisfies $\alpha(G)=\chi(G)$.}

   We have given a families of groups that are not DEG. A further problem is: \\
     {\bf Problem 4:}\ {\it Give more examples  that are not DEG. In particular, for any two positive integers $2\leq k<l$, give a group $G$ such that $\alpha(G)=k$ and $\gamma(G)=l$. }

  We    have proved that all finite simple groups and all finite nilpotent groups are DEG-groups.\\
    {\bf Problem 5:}\ {\it Give more families of groups that are DEG-groups and study the   properties of DEG-groups.}

An example is given in Section 4 to show  that  a minimum determining set in a DEG-group need not   be a minimum generating set. A natural problem is: \\
 {\bf Problem 6:}\ {\it Under what conditions, can all minimum determining sets of a  DEG-group   turn out to be minimum generating sets? }
\vskip 2mm
\noindent{\Large\bf Acknowledgments}
\vskip 2mm
 The first author would like to express his deep gratitude to Prof.
 Qinhai Zhang, Prof. Wenbin Guo, Prof. Xiuyun Guo, Prof. Xianhua Li,  Prof. Heguo Liu and Prof. Shenglin Zhou for their valuable comments, which   improve the quality of this paper a lot!

\end{document}